\documentclass[12pt]{article}
\usepackage{epsf}

\usepackage{amsmath,amsfonts,amssymb,amsthm,graphicx,color} \parskip1pt

\def\A{\= A}
\def\a{\= a}

\def\u{\= u}

\def\suls{{\'Sulvas\u tras}}

\def\sul{{\'Sulvas\u tra}}

\def\bsl{{Baudh\a yana}}
\def\asl{{\A pastamba}}
\def\ksl{{K\a ty\a yana}}
\def\msl{{M\a nava}}

\begin{document}  \date{}
\centerline{\large SOME CONSTRUCTIONS IN THE M\=ANAVA SULVAS\=UTRA\footnote{Text version of a talk given at the Annual Conference on History and Development of Mathematics - 2018, organized by the Sanskrit Academy, Madras in collaboration with and at Sri Chandrasekharendra Saraswati Vishwa Mahavidyalaya, Enathur, Kanchipuram, during November 27-29, 2018, under the auspices of the Indian Society for History of Mathematics.}} 

\bigskip
\centerline{S.G. Dani}
%\smallskip
\begin{center} 
{Distinguished Professor, UM-DAE Centre for Excellence in Basic Sciences,\\[1mm] Kalina, Santracruz,  Mumbai 400098}
%\maketitle
\end{center}
%{\large 
\vskip7mm
\noindent ---------------------------------------------------------------------------------------------------------
\noindent {\bf Abstract}: The \msl\ \sul, while less sophisticated than the other \suls, is seen to contain some mathematical ideas and constructions not found in the other \suls. Here we discuss some of these constructions and discuss their significance in the overall context of the \sul\ literature. 

\noindent 
---------------------------------------------------------------------------------------------------------

\smallskip
\noindent {\bf 1.  Introduction:} Among the works from the Vedic period that have come down to us, the \suls\ constitute a major source enabling understanding of that time concerning the mathematical aspects.   \suls\ were composed in aid of the activity around construction of {\it agnis} and {\it vedis} 
(fireplaces and altars) for performance of the {\it yajnas} which, it is needless to add here, had a very important role in the life of the Vedic people. The Vedic community was fairly heterogeneous, though with a shared tradition and body of knowledge, and there would have been numerous \suls, used by various local communities. Not surprisingly very few have survived. Of the handful of extant \suls, four are found to be significant from a mathematical point of view: \bsl\ \sul, \asl\ 
\sul, \msl\, \sul, \ksl\ \sul. 

While there is considerable
uncertainty about the time when the  \suls\ were composed, it has now become customary among the commentators to assign to their composition
 the period 800~-~200~BCE, with  \bsl\ \sul, believed to be the earliest, to be from around 800~-~500~BCE. It is also concluded from various
considerations that  \msl\ \sul\ is 
from a later period than \bsl\ \sul\ but is a little older than \asl\ \sul\
and considerably so compared to \ksl\ \sul; the ranges assigned typically are 650 - 300 BCE for \msl\ and \asl\ \suls\ and 400 - 200 for \ksl\ \sul. Despite being the oldest 
 \bsl\ \sul\ is found to be better organized and more elegant in its
presentation among all the four, while \msl\ \sul\ is least appealing
from these considerations. It  has also been the one to have 
received least attention in terms of editions, commentaries
etc., whether in traditional or in modern context, perhaps due to its lack of appeal.  The first modern
edition with English translation, due to Jeanette van Gelder
\cite{vG}, is only a little over 50 years old, while for the others similar
activity was undertaken well over a hundred years ago, in the nineteenth and early twentieth centuries. 

Notwithstanding its lack of appeal, there are some very interesting original observations 
in \msl\ \sul\ in terms of the mathematical content, which in the overall 
context seem to have not received adequate attention. I may also put in a  comment here that there seems to be a tendency among the scholars in the area to view the \suls\ body of knowledge mostly as a totality and  the special features of the individual \suls\ are scarcely highlighted, except at a superficial level, while,  on the other hand, there is no doubt that comparative studies between the individual \suls\ could throw a good deal of light on various aspects of the Vedic civilization, especially as the \suls\ are from different periods, and very likely also from different geographical regions in India. 
It is the aim
of this article to highlight some of the unique features of \msl\ \sul\
compared to the other \suls. 

\medskip
\noindent{\bf Circumferance of the circle}

\smallskip
During the ancient period, around the world the ratio of the circumference to the diameter
of the circle was thought to be 3,\footnote{One may wonder why the value for the ratio was taken to be 3 across
various cultures.  My hypothesis on the issue is that the idea of the
ratio being 3  
dates back to the time when humans were yet to think in terms of fractions (except 
perhaps for ``half'', which may have meant a substantial portion that is not nearly the whole - as commonly used even now in informal conversations -, 
rather than its precise value); 
it may be noted that while encounter with the circle, in the context
of wheels, is at least 
over 5000 years old, fractions seem to have appeared on the scene in a
serious way, 
in Indian as well as Egyptian cultures, only around the first millennium
BCE.  The ratio is thus 3 in the sense that it is
not 2 or 4, or even ``three hand half''.
The ingrained notion could have  developed into a belief (often tagged
also to religious authority). It was then not reconsidered for a
long time, even after fractions became part of human thought process. 
The episodes such as discussed here mark a departure from the past. 
} and the belief is also reflected in
one of the s\u tras in \bsl\ \sul; at one point there is an incidental reference to this, where a
circular pit ``with diameter 1 {\it pada} and circumference 3~{\it
  padas}'' is mentioned, indicating that the circumference was taken to
be 3 times the diameter.   The issue does not feature elsewhere in \bsl\ \sul\ and  in 
\asl\ and \ksl\ \suls. In the \msl\ \sul\ however one sees a
recognition that the assumption  is not correct. A verse  in  \msl\  (10.2.3.13 as per \cite{K} and 11.13 as per \cite{SB} numberings) states

\begin{quote}
{\tt vi\d skambha\d h pa${\tilde {\tt n}}$cabh\a ga\'sca vi\d skambhastrigu\d na\'sca ya\d h}$|$

{\tt sa ma\d n\d dalaparik\d sepo na v\a lamatiricyate}$||$

\medskip
``a fifth of the diameter and thrice the diameter is the circumference
of a circle, not a hair-breadth remains.''

\end{quote}
``vi\d skambha", which also means supporting beam or bolt or bar of a door (see \cite{MW}, \cite{A}), was the technical term used for the diameter of a circle. ``ma\d ndala" stands for the circle and ``parik\d sepa\d h" is the term for the circumference.  
Even though the value described is considerably off the mark, the fact
of recognition of the ratio being strictly greater than 3 is worth taking
note of, and so is the apparent exultation over the finding.\footnote{In \cite{vG}, and in \cite{K} following it, the verse is wrongly interpreted as concerning determination of a square with the same area as the given circle: the translation of the verse is given as ``Dividing the diameter of the circle into five parts and then individual parts into three parts each (thus dividing the diameter into 15 parts  and taking away two parts) yields the side of a square with the same area as the circle. This is accurate to a hair-breadth." If the translation in the first part were to be correct then it would correspond to the formula for the side of the square with area equal to that of a given circle is $\frac {13}{15}$ seen in \bsl\ (at 1.60) and also in \asl\ and \ksl\ \suls. The translation however is quite erroneous in many respects: occurrence of the word ``vi\d skambha" twice readily shows that it is not the individual parts that are being subdivided, and there is no reference at all to taking away two parts from the 15 subdivided parts. Besides, ``parik\d shepa" unambiguously corresponds to circumference, with the verb ``parikship" meaning ``to surround", ``to encircle" etc. (see \cite{MW}, \cite{A}), and not the area. It appears that having difficulties in interpreting the verse the translator chose to relate it to the $\frac {13}{15}$ formula seen in the other \suls. The translation in \cite{SB} on the other hand is along the lines described here.} 

It may be recalled here that in the Jaina tradition a similar recognition is seen in {\it s\u
  ryaprajn\a pti} (believed to be from 5 th century BCE), where the
classical value 3 for the ratio is 
recalled and discarded in favour of another value $\sqrt {10}$. The values could thus be contemperanous, but evidently unrelated from a historical point of view, especially on account of the substantial difference in the values proposed, in numerical as well as structural terms.  

A brief description of the location of the verse in the body of the \msl\ \sul\  would be in order here, to place the verse in context. 
Section 10.3 in which the verse occurs, at 10.3.2.13,\footnote{Actually the 10 is superfluous in these numbers, since the whole of \msl\ \sul\ is covered in sections of Chapter 10; the numbering has to do with the translation of {\it M\=anava\ \'Srautas\u tras} in \cite{vG}, in which the \sul\ appears as Chapter~10.}  is the last of the three sections in \msl\ \sul, referred by the sulvak\=ara as {\it vai\d s\d nava}; the significance of the name, and association with {\it Vi\d s\d nu}, if any,  is not clear from the contents of the section. The general narrative in the part containing the verse concerns description of construction of  {\it vedis}.   Interestingly, after talking about the volume of the vedi called {\it \'samitra vedi} the s\=utrak\=ara states 

\begin{quote}
{\tt \=ay\=amam\=ay\=amagu\d nam vist\=ara\d m vistare\d na tu~$|$

samasya vargm\=ulam yat tatkar\d na\d m tadvido vidu\d h $|$}

\medskip
``Multiply the length by the length and the width by the width. It is known that  adding them and taking the square root gives the hypotenuse.'' 
\end{quote}

The reader would recognize this statement as an equivalent form of what is called the Pythagoras theorem (!), with the figure in question (not specified in the verse) being the rectangle.\footnote{\cite{K} also mentions the right-angled {\it triangle} in this respect but there is no evidence, on the whole, of the \suls\ discussing right angled triangles.}
It may also be noted that the statement is in quite a different form than in the other \suls; in a way, while the other \suls\ seem to be referring to geometric principle involved, considering in particular  the areas of the squares over the respective sides, the exposition here is seen to be focussed on {\it computing} the the size of the hypotenuse from the sizes of the sides, without specific reference to the underlying geometry. 

A few verses down from there, which concern practical details about  the vedis and the performance of yajna, we are led to another important mathematical statement, involving now the construction of a circle with the same area as a square.\footnote{It is argued in \cite{H} that 10.3.2.10 gives rules both for squaring the circle, and circling the square, with the latter being the same as Baudhayana's rule discussed earlier. The rule for the other direction, according to the interpretation in \cite{H} is that given a circle, the perpendicular bisector of the equilateral triangle with the diameter of the circle as the side, is the side of a square with the same area as the circle. The argument involves an emendation of the extant text, which the author justifies also on considerations of grammar, but with it many aspects which are unclear from the earlier translations from \cite{vG} and \cite{SB} become clearer. As noted in \cite{H} the above mentioned rule for squaring the circle is unique to \msl\ \sul. The rule however is not very accurate.} In the vedic literature this issue concerns constructing the {\it \=ahavan\=iya}, which is a square and {\it g\=arhapatya} which is circular with the same area\footnote{In \cite{D} concerning the motivation for considering the problem of circling the square, I had made a reference to the {\it rathacakraciti}, which however seems to be an inadequate explanation - the primary motivation for the problem is very likely to have been the equality of the areas of {\it \=ahavan\=iya} and {\it g\=arhapatya}. }; along with there is also the semicircular figure with the same area to be constructed for the {\it dak\'si\d n\=agni}. The method described here for finding a circle with the same area as a given square is the same as given in \bsl\ \sul\ in geometrical content, but formulated with a difference: in an isosceles triangle produced by the diagonals of the square, extend the perpendicular (as much as the semidiagonal side of the triangle) and of the extra part of the semidiagonal (beyond the side) adjoin a third part of it to the part within the square, to get the radius of the circle.   As is well-known (see in particular \cite{D} for a discussion on this) this is not very accurate, but is interesting as an approximate construction.  This is followed by two  verses which concern doubling of area when measure of a side is replaced by that of the diagonal of the square. This is evidently related in this context with the construction of the {\it dak\d si\d n\=agni}, though it has not been explicitly mentioned, and has also not been brought out in the translations in \cite{K} and \cite{SB}. 

And then comes the cited verse for the circumference of the circle! What is the relevance that we can identify? We see that some circles have appeared on the scene, though what is involved about them are the {\it areas}. Nothing in the context warrants, apparently, consideration of the circumference. However, having got to the circles, seems to have inspired the author to mention, and that too with some gusto, something interesting that he {\it had realized}, namely that the circumference is not just 3 times the diameter as people thought, but more than that, and one would have a safe estimate by adding one-fifth of the diameter.  Thus the statement (like much else actually, but
 it bears emphasis here) appears to be side input, 
from which it would be difficult to draw further inference about the thought process that may be involved.

Indeed, one may wonder why the sulvak\=ara chose the value $3\frac 15$ for the correction, rather than something that would have been better, specifically like $\frac 16$, if not $\frac 17$.  From the context, and the value itself, it is clearly an ad hoc value being adopted, essentially in the context of becoming aware of the classical value of $3$ for the ratio is not satisfactory, and that {\it something remains}. I may reiterate here in this respect that the verse notes that ``not a hair-breadth {\it remains}", which is what ``{\it atiricyate}" corresponds to, with the verb ``{\it atiric}" meaning ``to be left with a surplus" (see \cite{MW}), and is strictly not a reference accuracy in terms of both lower and upper estimates (as treated, for instance, in \cite{K}). But having recognized that the value should be more than 3, why and how did $3\frac 15$ come to be chosen for it. The value $3\frac 17$, which would be appropriate in hindsight,  would perhaps would have been rather odd (lacking in aesthetic appeal, which is often a consideration while making ad hoc choices) to think about at that time. However, why not, say $3\frac 16$, which would have been much closer to the correct value? Thinking of a sixth would seem simpler and natural compared a fifth part, it being half of the a third, and division into three parts is easier operationally, than into five parts, and then halving would of course be the trivial next step. 

The s\=utrakara however prefers to consider division into five parts. A clue for this seems to lie in the decimal place value system of representation of numbers (writing numbers to base 10, as we do now). For a number written in this system it is much easier to compute its fifth part than the third, or any other, part. Indeed, \msl\ \sul\ shows preference to using decimally convenient divisions in other contexts as well.  The verse following the cited one, for the circumference describes the size of a square inscribed in a circle, viz. with vertices on the circumference. It may be noted that the desired size would be $1/\sqrt 2$ times the diameter of the circle. The prescription given is to divide the diameter into 10 parts and take away 3 parts; thus $7/10$ is used as as a(n approximate) value for $1/\sqrt 2$. Actually for $\sqrt 2$ there was a standard approximate value $17/12$ adopted in the  \suls, according to which the desired ratio would be 12 out of 17 parts, which would be more accurate, but \msl\ adopts the proportion 7 out of 10, suggesting preference for decimal division. In the verse for a new construction for circling the square, which we shall discuss in the next section, there is a division into 5 parts involved. It may also be recalled here that the major large unit involved in the \suls\ is {\it puru\d sa} and there is a subunit {\it aratni}, which is a 5th part of {\it puru\d sa}. This may also be looked upon as a factor, related to the use of the decimal system, which would have encouraged considering division into 5 parts.   I may also recall here that in   the construction of various vedis that are described, division by 5 is involved in many computations.   

We conclude this discussion with another small related observation. Granting that the value of the circumference to diameter ratio was recognized by the sulvak\=ara as being greater than 3, and that he looked for additional decimal parts after which ``nothing will remain", $1/5$ is the right choice; $1/10$ would have been closer, but it is {\it less} than the correct value.    

\medskip
\noindent {\bf 3. Circling the Square}

\smallskip
As noted in the last section the problem of circling the square, viz. of describing a circle with the same area as a given square, had attained considerable importance in the \suls\ period. It may be emphasized that the framework envisaged for the problem is quite different from the analogous problems in Greek mathematics, where the constructions were sought to be performed with {\it only} the ruler and compass, and any comparison of the achievements of the ancient Indians, in the context of the Greeks ``not having been successful" with the problem are facile and irrelevant. The constructions given are important in terms of historical development of mathematical ideas, and need to be viewed only as such. 

We have gone over the geometric construction given in \bsl\ \sul\ for drawing a circle with the same as a given square. As noted there, the result it produces is not very accurate, and in fact involves an error of the order of 1.7 percent (see \cite{D} for more details in this respect). 
In course of time the suspicions would have 
gained weight, serving as motivation to look for
an alternative construction, and one seems to  find such an attempted construction in  
\msl\ \sul, which we shall now discuss.

The construction in question is described in a verse which follows right after the contents discussed in the last section here (in  10.3.2.15 as per the numbering of  \cite{K} and 11.15 of \cite{SB}). The verse is 

\begin{quote}
{\tt caturasra\d m navadh\=a kury\=at dhanu\d h ko\d tistridh\=atridh\=a$|$\\
utsedh\=atpa$\tilde{\tt n}$cama\d m lumpetpur\={\i}\d seneh t\=avatsama\d m $||$}

\end{quote}

The first part of the verse may be translated, quite unambiguously, as

\begin{quote}
  Divide the square into nine parts, (by) dividing the horizontal and vertical sides into three parts each.
 \end{quote}
 
Unfortunately, arriving at the right translation of the rather terse second part of the verse, and its interpretation, call for additional  inputs of contextual nature, and want of these seems to have confused earlier translators: in  \cite{SB} the authors translate the second part 
as 

\begin{quote}
``drop out the fifth portion (in the centre) and fill it up with loose earth."

\end{quote}
and in the commentary section they comment 

\begin{quote}
``Possibly these are  not problems 
  of quadrature of the circle. Ordinary squares are drawn without any
  mathematical significance.'' 
\end{quote}
 The comment seems quite unwarranted, though it may be emphasized there that nothing in the verse specifically indicates that it does concern a quadrature formula, or procedure towards one.  In  \cite{K}, following \cite{vG}, the second part is interpreted (in Marathi and Hindi equivalents) as

\begin{quote}  
  ``from the part jutting out take {away} one-fifth part and draw a circle with the remaining part as the radius". 
\end{quote}  
 Here the word ``{\it utsedha}" is interpreted to mean the part of the  trisectors (arrived at in the first part) that is jutting out on either side of the square, until meeting the circle passing through the vertices of the square. One fifth of that is subtracted from the segment of the trisector upto its midpoint, and the remaining part is taken as the radius of the prescribed circle. Implementing the procedure accordingly, they calculate the radius; it however turns out to be much too large for the corresponding circle to have the same area as the square, thus putting the interpretation into question, but the matter is left  at that, with no comment.  
  
Another interpretation of the verse was given by R.C. Gupta \cite{G} (see also \cite{G2}). Here ``utsedha" is interpreted to mean ``height" and is associated with the ``height", viz. the radius, of the semicircle from the circle through the vertices of the square. Thus the author infers that the radius of the prescribed circle is meant to be $\frac 45$th of the circumscribing circle. With this interpretation the area of the circle produced, starting with the unit square, works out to be $8\pi/25$, and thus the procedure corresponds to a value of $\pi$ as $25/8$. This is a good value by the \suls\ standards. However the interpretation is unsatisfactory in various ways: First and foremost, the interpretation does not involve the first part of the verse at all. It is inconceivable that the sulvak\=ara first asks you to elaborately divide the square into nine parts, and in the following line gives a procedure for the quadrature problem which has nothing to with the subdivision.  Secondly, in the second part if it was just the radius of the circumscribing circle to be used as a reference why would it be referred to with the unusual word ``utsedha", which does not occur anywhere else in \msl\ (or in other \suls), rather than in terms of the diameter of the circle, which is something that occurs so frequently in \suls\ geometry. 

For a faithful interpretation of the verse it seems imperative that it must involve the trisectors of square introduced in the first part; also {\it utsedha} must have something to do with the trisectors and the choice of the unusual term must have to do with that the trisectors also do not occur anywhere else.  Thus it would seem that the interpretation in \cite{vG} and \cite{K} is on the right track inasmuch as it focusses on considering individually the lines trisecting the given square along each of the sides,  extended upto the circle passing through the vertices of the circle. The circle is indeed being described in terms of certain points on these lines. The main difficulty however seems to be in understanding which points are meant. Evidently the interpretation with regard to the points, and how they are to be used (see more on this below), adopted in \cite{vG} and \cite{K} does not seem to the right one, as it is way off the mark.  

The overall formulations and symmetry considerations suggests that we are to pick two points on each of the four lines that trisect the square along a side, located symmetrically (and hence at the same distance from the centre of the square) and the circle through these points is the desired circle; this in a way explains the explication through ``covering with loose earth", as the totality of the eight points is indicative of a circle which is what is to be covered. Now, which are the two points on each of the lines? One would be in a better position to figure out what the sulvak\=ara's line of thought, if one keeps in mind the \bsl\ construction of the circle, described earlier. Recall that there the bisector of the square is extended until meeting the circle through the vertices, and $\frac 13$rd of the part is added to the segment within, to get the radius of the circle that is sought after; one can alternatively think of this as identifying the point through which the circle should pass (the centre of the circle is of course understood to be the centre of the square). The new idea now is that instead of the bisectors of the squares we are considering  the trisectors. On the bisectors the points in question were chosen to be at $\frac 13$rd of the jutting out part, {\it from the side of the square}. One now needs to look for a similar number, and the analogous point on the trisectors, to complete the analogous construction. The number is picked to be $\frac 15$th; the choice could have been based on intuition, and the point is now meant to be on the trisector at $\frac 15$th of the jutting out part, from the side of the square. At this point the analogy with the \bsl\ construction throws open, to our minds, two possibilities, one is to take the segment of the trisector upto the point thus constructed either from the midpoint of the trisector, or the centre of the square; in the case of the \bsl\ construction with the bisector the two coincide, but here they are different. For some reason in \cite{K} the former interpretation is favoured (with respect to the point picked there, on which was commented upon above). However, viewed in the full context, it is the other interpretation that may be seen to be more appropriate. The sulvak\=aras do not in general try describe a number for the radius, but a region to be covered determined by some point (or a collection of points), and secondly in the overall context of the description of the construction the midpoint of the trisector has no relevance (and has not been referred to). Once these points are noted, the inference would be that  the prescribed circle passing through the point(s) as above on the trisectors, at one-fifth of the jutting out part from the side of the square. One may now rewrite the interpretation of the second part,  referring to the collection of the 8 points, as 

\begin{quote}

``on the parts jutting out mark the points at one-fifth (from the square) and draw the circle through them."

\end{quote}
    
A simple calculation shows that for a square with 
unit side length  the radius of the circle is 
$\frac 12 \{1+\frac 15 (\frac {\sqrt {17}}3 -1)\}^2+\frac 1{18}$, and this yields
the area of the circle to be $0.994...$, a much more accurate value
compared to the earlier one, with an error of only about $0.5\%$, in place of $1.7\%$ (see \cite{D} for details of the calculations and other related comments). Thus from a mathematical point of view this turns out to be a good choice. We see also that it emerges naturally as a  generalization of the \bsl\ construction in terms of development of ideas. Thus it seems reasonable to expect that this is what the sulvak\=ara had in mind. The interpretation incorporates all the components of the verse, and all the ingredients needed in the formulation may be seen to be present in the verse, in their natural order. The author is hopeful that the interpretation would be confirmed to be valid by expert Sanskritists, from a linguistic point of view, possibly after some emendation that could be could be justified based on considerations of corruption on account of one or other factors.

\bigskip
\noindent 
S.G. Dani\\
UM-DAE Centre for Excellence in Basic Sciences\\
Campus of University of Mumbai, Kalina, Santrcruz\\                  
Mumbai 400 098\\
India 

\noindent E-mail: {\tt shrigodani@gmail.com}

\end{document}